\RequirePackage{ifpdf}
\ifpdf 
\documentclass[pdftex]{sigma}
\else
\documentclass{sigma}
\fi

\usepackage[mathscr]{eucal}

\begin{document}

\allowdisplaybreaks

\renewcommand{\thefootnote}{$\star$}

\renewcommand{\PaperNumber}{078}

\FirstPageHeading

\ShortArticleName{Harmonic Analysis on Quantum Complex Hyperbolic Spaces}

\ArticleName{Harmonic Analysis\\ on Quantum Complex Hyperbolic Spaces\footnote{This paper is a
contribution to the Special Issue ``Relationship of Orthogonal Polynomials and Special Functions with Quantum Groups and Integrable Systems''. The
full collection is available at
\href{http://www.emis.de/journals/SIGMA/OPSF.html}{http://www.emis.de/journals/SIGMA/OPSF.html}}}

\Author{Olga BERSHTEIN and Yevgen KOLISNYK}

\AuthorNameForHeading{O.~Bershtein and Ye.~Kolisnyk}

\Address{Institute for Low Temperature Physics and Engineering,\\ 47 Lenin Ave., 61103, Kharkov, Ukraine}
\Email{\href{mailto:bershtein@ilt.kharkov.ua}{bershtein@ilt.kharkov.ua}, \href{mailto:evgen.kolesnik@gmail.com}{evgen.kolesnik@gmail.com}}

\ArticleDates{Received April 30, 2011, in f\/inal form August 10, 2011;  Published online August 18, 2011}

\Abstract{In this paper we obtain some results of harmonic analysis on quantum complex hyperbolic spaces. We
introduce a quantum analog for the Laplace--Beltrami operator and its radial part. The latter appear to be
second order $q$-dif\/ference operator, whose eigenfunctions are related to the Al-Salam--Chihara polynomials.
We prove a Plancherel type theorem for~it.}


\Classification{17B37; 20G42; 81R50; 33D45; 42C10}

\section{Introduction}

Consider the group $SU_{n,m}$ and its homogeneous space
$\mathscr{H}_{n,m}=SU_{n,m}/S(U_{n,m-1}\times U_1)$. The latter is called a
complex hyperbolic space. The Faraut paper \cite{Faraut} on such
pseudo-Hermitian symmetric spaces has a great impact into the theory of
semisimple symmetric spaces of rank~1. Also there are numerous papers of
Molchanov, van Dijk and others (see \cite{Molch1,Molch2,Dijk-Sh} and
references therein) on representation theory related to these symmetric
spaces and harmonic analysis on them. In particular, there is the celebrated
Penrose transform which enables to relate classical bounded symmetric
domains and complex hyperbolic spaces.

In this paper we develop harmonic analysis on quantum complex hyperbolic spaces. Recall that the related polynomial
algebras were introduced in the early paper by Faddeev, Reshetikhin, and Takhtadjan~\cite{RTF}. Unfortunately,
there was no further inquiry.

Nearly 10 years ago L.~Vaksman and his team started the theory of quantum bounded symmetric domains. Their
approach enables to formulate and solve various problems on noncommutative complex and harmonic analysis in
these domains, geometric realizations of representations of quantum groups~\cite{V}. Their paper~\cite{Penrose} establishes the missed link between quantum bounded symmetric domains and introduces this
project.

The initial notions of function theory on quantum complex hyperbolic spaces $\mathscr{H}_{n,m}$ were
introduced in~\cite{BS}. Namely, a quantum analog $\mathscr{D}(\mathscr{H}_{n,m})_{q,\mathfrak{k}}$ for the
algebra of $U \mathfrak{k} = U \mathfrak{s(gl_n \times gl_m)}$-f\/inite smooth functions on $\mathscr{H}_{n,m}$ with compact support
and an $U_q \mathfrak{su}_{n,m}$-invariant integral $\int d\nu_q$ on it were constructed. We recall basic
notations from the quantum group theory and \cite{BS} in Sections~\ref{sec_1},~\ref{sec_2},~\ref{finite}.
Section~\ref{sec_4} is devoted to the subalgebra of $U_q \mathfrak{k}$-invariant f\/inite functions.

In Sections~\ref{sec_5} and~\ref{sec_6} we introduce a quantum analog $\square$ for the Laplace--Beltrami operator on~$\mathscr{H}_{n,m}$. Also we consider its radial part $\square^{(0)}$, i.e., its restriction to the space
$L^2\big(d\nu_q^{(0)}\big)$, where $\int \cdot d\nu_q^{(0)}$ is the restriction of $\int d\nu_q$ to the space of $U_q
\mathfrak{k}$-invariant elements of $\mathscr{D}(\mathscr{H}_{n,m})_{q,\mathfrak{k}}$. The latter operator
naturally appears to be $q$-dif\/ference operator and is related to a three-diagonal Jacobi matrix. Section
\ref{sec_7} is devoted to generalized eigenfunctions of $\square^{(0)}$ and lead us towards an expected yet
remarkable appearance of Al-Salam--Chihara polynomials. A spectral theorem for~$\square^{(0)}$ in Section~\ref{sec_9} is obtained as a~corollary to well known results on these polynomials (see~\cite{koekoek}).

\section{Preliminaries on quantum group theory}\label{sec_1}

Everywhere in the sequel we suppose $q\in(0,1)$. All algebras are associative and unital.

The Hopf algebra $U_q\mathfrak{sl}_{N}$ is given by its generators $K_i$, $K_i^{-1}$, $E_i$, $F_i$,
$i=1,2,\ldots,N-1$, and the relations:
\begin{gather*}
K_iK_j=K_jK_i,\qquad K_iK_i^{-1}=K_i^{-1}K_i=1,
\\
K_iE_i=q^2E_iK_i,\qquad K_iF_i=q^{-2}F_iK_i,
\\
K_iE_j=q^{-1}E_jK_i,\qquad K_iF_j=q F_jK_i, \qquad |i-j|=1,
\\
E_iF_j-F_jE_i=\delta_{ij}\frac{K_i-K_i^{-1}}{q-q^{-1}},
\\
E_i^2E_j-\big(q+q^{-1}\big)E_iE_jE_i+E_jE_i^2=0,\qquad|i-j|=1,
\\
F_i^2F_j-\big(q+q^{-1}\big)F_iF_jF_i+F_jF_i^2=0,\qquad|i-j|=1,
\\
[E_i,E_j]=[F_i,F_j]=0,\qquad|i-j|\ne 1.
\end{gather*}
 The comultiplication $\Delta$, the antipode $S$, and the counit
$\varepsilon$ are def\/ined on the generators by
\begin{gather*}
\Delta(E_i)=E_i\otimes 1+K_i\otimes E_i,\qquad \Delta(F_i)=F_i\otimes K_i^{-1}+1\otimes F_i,\qquad
\Delta(K_i)=K_i\otimes K_i,
\\
S(E_i)=-K_i^{-1}E_i,\qquad S(F_i)=-F_iK_i,\qquad S(K_i)=K_i^{-1},
\\
\varepsilon(E_i)=\varepsilon(F_i)=0,\qquad\varepsilon(K_i)=1,
\end{gather*} see \cite[Chapter 4]{Jant}.

We need also the Hopf algebra $\mathbb{C}[SL_{N}]_q$ of matrix elements of f\/inite dimensional weight
$U_q\mathfrak{sl}_{N}$-modules. Recall that $\mathbb{C}[SL_{N}]_q$ can be def\/ined by the generators $t_{ij}$,
$i,j=1,\dots,N$ (the matrix elements of the vector representation in a weight basis) and the relations
\begin{alignat*}{3}
&t_{ij'}t_{ij''}=qt_{ij''}t_{ij'},\qquad && j'<j'',&
\\ &t_{i'j}t_{i''j}=qt_{i''j}t_{i'j},\qquad && i'<i'',&
\\ &t_{ij}t_{i'j'}=t_{i'j'}t_{ij},\qquad && i<i'\;\&\;j>j', &
\\ &t_{ij}t_{i'j'}=t_{i'j'}t_{ij}+\big(q-q^{-1}\big)t_{ij'}t_{i'j},\qquad
&& i<i'\;\&\; j<j',&
\end{alignat*}
together with one more relation
\[
\det\nolimits_q\mathbf{t}=1,
\]
where $\det\nolimits_q\mathbf{t}$ is a $q$-determinant of the matrix $\mathbf{t}=(t_{ij})_{i,j=1,\dots,N}$:
\[
\det\nolimits_q\mathbf{t}=\sum\limits_{s\in S_{N}}(-q)^{l(s)}t_{1 s(1)} t_{2s(2)}\cdots t_{Ns(N)},
\]
with $l(s)=\mathrm{card}\{(i,j)|i<j\;\&\;s(i)>s(j)\}$. The algebra $\mathbb C[SL_N]_q$ is endowed with the
standard structure of $U_q^{\rm op}\mathfrak{sl}_N \otimes U_q \mathfrak{sl}_N$-module algebra (here 'op'
ref\/lects the fact that we should change the multiplication by the opposite one).

Let also $U_q\mathfrak{su}_{n,m}$, $m+n=N$, denotes the Hopf $*$-algebra $(U_q \mathfrak{sl}_N,*)$ given by
\[
(K_j^{\pm 1})^*=K_j^{\pm 1},\qquad E_j^*=
\begin{cases}
K_jF_j,& j \ne n,
\\ -K_jF_j,& j=n,
\end{cases}\qquad F_j^*=
\begin{cases}
E_jK_j^{-1},& j \ne n,
\\ -E_jK_j^{-1},& j=n,
\end{cases}
\]
with $j=1,\ldots,N-1$ \cite{RTF, polmat}.

Recall the notion of an algebra of `regular functions on the quantum principal homogeneous space' $X$
constructed in \cite{polmat}. Put $\operatorname{Pol}(\widetilde{X})_q\overset{\mathrm{def}}{=}
(\mathbb{C}[SL_N]_q,*)$, where the involution $*$ is def\/ined by
\begin{gather}\label{def_inv}
t_{ij}^*=\mathrm{sign}[(i-m-1/2)(n-j+1/2)](-q)^{j-i}\det\nolimits_qT_{ij}.
\end{gather}
Here $\det_q$ is the quantum determinant \cite{Ch-P}, and the matrix $T_{ij}$ is derived from the matrix
$\mathbf{t}=(t_{kl})$ by discarding its $i$'s row and $j$'s column. In \cite{polmat} it is proved that
$\operatorname{Pol}(\widetilde{X})_q$ is a $U_{q}\mathfrak{su}_{n,m}$-module $*$-algebra.

\section[A $*$-algebra $\operatorname{Pol}(\mathscr{H}_{n,m})_{q}$]{A $\boldsymbol{*}$-algebra $\boldsymbol{\operatorname{Pol}(\mathscr{H}_{n,m})_{q}}$}\label{sec_2}

Let $m,n\in\mathbb{N}$, $m\ge 2$, and $N\overset{\mathrm{def}}{=}n+m$. Recall that the classical complex
hyperbolic space $\mathscr{H}_{n,m}$ can be obtained by projectivization of the domain
\[
\widehat{\mathscr{H}}_{n,m}=\left\{(t_1,\ldots,t_N)\in\mathbb{C}^N\; \Big|\;
-\sum^n_{j=1}|t_j|^2+\sum^N_{j=n+1}|t_{j}|^2>0 \right\}.
\]

Now we pass from the classical case $q=1$ to the quantum case $0<q<1$. Let us consider the well known
\cite{RTF} $q$-analog of the polynomial algebras. Let $\operatorname{Pol}(\widehat{\mathscr{H}}_{n,m})_q$
be the unital $*$-algebra with the generators $t_1,t_{2},\ldots,t_N$ and the commutation relations as
follows:
\begin{gather*}
t_it_j = qt_jt_i,\qquad i<j,
\\ t_it_j^* = qt_j^*t_i,\qquad i\ne j,
\\ t_it_i^* = t_i^*t_i+\big(q^{-2}-1\big)\sum_{k=i+1}^Nt_kt_k^*,\qquad i>n,
\\ t_it_i^* = t_i^*t_i+\big(q^{-2}-1\big)\sum_{k=i+1}^nt_kt_k^*-
\big(q^{-2}-1\big)\sum_{k=n+1}^Nt_kt_k^*,\qquad i\le n.
\end{gather*}
Obviously, \[c=-\sum_{j=1}^nt_jt_j^*+\sum_{j=n+1}^Nt_jt_j^*\] is central in
$\operatorname{Pol}(\widehat{\mathscr{H}}_{n,m})_q$. Moreover, $c$ is not a zero divisor in
$\operatorname{Pol}(\widehat{\mathscr{H}}_{n,m})_q$. This allows one to embed the $*$-algebra
$\operatorname{Pol}(\widehat{\mathscr{H}}_{n,m})_q$ into its localization
$\operatorname{Pol}(\widehat{\mathscr{H}}_{n,m})_{q,c}$ with respect to the multiplicative system
$c^{\mathbb{N}}$.

The $*$-algebra $\operatorname{Pol}(\widehat{\mathscr{H}}_{n,m})_{q,c}$ admits the following bigrading:
\[\deg t_{j}=(1,0),\qquad\deg t_j^*=(0,1),\qquad j=1,2,\ldots,N.\] Introduce the notation
\[
\operatorname{Pol}(\mathscr{H}_{n,m})_q=\left\{ f\in
\operatorname{Pol}(\widehat{\mathscr{H}}_{n,m})_{q,c} \; \big|\; \deg f=(0,0)\right\}.
\]
This $*$-algebra $\operatorname{Pol}(\mathscr{H}_{n,m})_q$ will be called the algebra of regular functions on
the quantum hyperbolic space.

We are going to endow the $*$-algebra $\operatorname{Pol}(\mathscr{H}_{n,m})_q$ with a structure of
$U_{q}\mathfrak{su}_{n,m}$-module algebra~\cite{Ch-P}. For this purpose, we embed it into the
$U_{q}\mathfrak{su}_{n,m}$-module $*$-algebra $\operatorname{Pol}(\widetilde{X})_q$.

By a $q$-analog of the Laplace expansion of $\det_q \mathbf{t}$ along the f\/irst row \cite[Section~9.2]{Kl-Sch} and
\eqref{def_inv}, one can obtain from $\det_q \mathbf{t}=1$ that
\[
-\sum_{j=1}^nt_{1j}t_{1j}^*+\sum_{j=n+1}^Nt_{1j}t_{1j}^*=1.
\]
Thus the map $J:t_j\mapsto t_{1j}$, $j=1,2,\ldots,N$, admits a unique extension to a homomorphism of
$*$-algebras $J:\operatorname{Pol}(\widehat{\mathscr{H}}_{n,m})_{q,c}\to
\operatorname{Pol}(\widetilde{X})_q$. Its image will be denoted by
$\operatorname{Pol}(\widetilde{\mathscr{H}}_{n,m})_q$. It is easy to verify that the $*$-algebra
$\operatorname{Pol}(\mathscr{H}_{n,m})_q$ is {\it embedded} this way into
$\operatorname{Pol}(\widetilde{\mathscr{H}}_{n,m})_q$ and its image is just the subalgebra in
$\operatorname{Pol}(\widetilde{\mathscr{H}}_{n,m})_q$ generated by $t_{1j}t_{1k}^*$, $j,k=1,2,\ldots,N$
(recall that $c$ goes to $\det_q \mathbf{t}=1$). In what follows we will identify
$\operatorname{Pol}(\mathscr{H}_{n,m})_q$ with its image under the map $J$.

Consider the subalgebra $U_q \mathfrak{s(gl}_1 \times \mathfrak{gl}_{N-1})$ generated by $K_i^{\pm 1}$,
$i=1,\dots,N-1,$ $E_j$, $F_j$, $j =2,\dots,N-1.$ By obvious reasons,
\[
\operatorname{Pol}(\mathscr{H}_{n,m})_q=\left\{f \in \operatorname{Pol}(\widetilde{X})_q \;\big|\;
L(\xi)f=\varepsilon(\xi)f, \xi \in U_q \mathfrak{s(gl}_1 \times \mathfrak{gl}_{N-1})\right\},
\]
where $L$ is the left action of $U_q^{\rm op} \mathfrak{sl}_N$ in $\operatorname{Pol}(\widetilde{X})_q$.

Let $I_\varphi$, $\varphi\in\mathbb{R}/2\pi\mathbb{Z}$, be the $*$-automorphism of the $*$-algebra
$\operatorname{Pol}(\widetilde{\mathscr{H}}_{n,m})_q$ def\/ined on the generators $\{t_j\}_{j=1,\ldots,N}$ by
\begin{gather}\label{I_phi}
I_\varphi: \ t_j\mapsto e^{i\varphi}t_j.
\end{gather}
We
use the notation $t_j$ instead of $t_{1j}$ for the generators of
$\operatorname{Pol}(\widetilde{\mathscr{H}}_{n,m})_q$.

Then one more description of $\operatorname{Pol}(\mathscr{H}_{n,m})_q$ is as follows:
\[
\operatorname{Pol}(\widetilde{\mathscr{H}}_{n,m})_q = \left\{ f\in
\operatorname{Pol}(\widetilde{\mathscr{H}}_{n,m})_q\; \big| \;I_\varphi(f)=f\text{ \ for all \
}\varphi\right\}.
\]

At the end of this section we list explicit formulas for the action of $U_{q}\mathfrak{su}_{n,m}$ on
$\operatorname{Pol}(\widetilde{\mathscr{H}}_{n,m})$:
\begin{gather*}
E_jt_i =
\begin{cases}
q^{-1/2}t_{i-1}, & j+1=i,
\\ 0, & \text{otherwise},
\end{cases}\qquad
F_jt_i =
\begin{cases}
q^{1/2}t_{i+1}, & j=i,
\\ 0, & \text{otherwise},
\end{cases}
\\ K_j^{\pm 1}t_i  =
\begin{cases}
q^{\pm 1}t_i, & j=i,
\\ q^{\mp 1}t_i, & j+1=i,
\\ t_i, & \text{otherwise}.
\end{cases}
\end{gather*}

\section[Algebras of generalized and finite functions on the quantum $\mathscr{H}_{n,m}$]{Algebras of generalized and f\/inite functions\\ on the quantum $\boldsymbol{\mathscr{H}_{n,m}}$}\label{finite}

Let us construct a faithful $*$-representation $T$ of $\operatorname{Pol}(\mathscr{H}_{n,m})_q$ in a
pre-Hilbert space $\mathscr{H}$ (our method is well known; see, for example, \cite{polmat}).

The space $\mathscr{H}$ is a linear span of its orthonormal basis
$\{e(i_1,i_2,\ldots,i_{N-1})\, |\,i_1,\ldots,i_n\in-\mathbb{Z}_+$; $i_{n+1}, \ldots,i_{N-1}\in\mathbb{N}\}$.

The $*$-representation $T$ is a restriction to $\operatorname{Pol}(\mathscr{H}_{n,m})_q$ of the
$*$-representation of $\operatorname{Pol}(\widetilde{\mathscr{H}}_{n,m})_q$ def\/ined by
\begin{gather*}
T(t_j)e(i_1,\ldots,i_{N-1})  = q^{\sum\limits_{k=1}^{j-1}i_k} \big(q^{2(i_j-1)}-1\big)^{1/2}
e(i_1,\ldots,i_j-1,\ldots,i_{N-1}),
\\ T(t_j^*)e(i_1,\ldots,i_{N-1})  =
q^{\sum\limits_{k=1}^{j-1}i_k} \big(q^{2i_j}-1\big)^{1/2} e(i_1,\ldots,i_j+1,\ldots,i_{N-1}),
\end{gather*}
for $j\le n$,
\begin{gather*}
T(t_j)e(i_1,\ldots,i_{N-1})  = q^{\sum\limits_{k=1}^{j-1}i_k} \big(1-q^{2(i_j-1)}\big)^{1/2}
e(i_1,\ldots,i_j-1,\ldots,i_{N-1}),
\\ T(t_j^*)e(i_1,\ldots,i_{N-1})  =
q^{\sum\limits_{k=1}^{j-1}i_k} \big(1-q^{2i_j}\big)^{1/2} e(i_1,\ldots,i_j+1,\ldots,i_{N-1}),
\end{gather*}
for $n<j<N$, and, f\/inally,
\begin{gather*}
T(t_N)e(i_1,\ldots,i_{N-1})  = q^{\sum\limits_{k=1}^{N-1}i_k}e(i_1,\ldots,i_{N-1}),
\\ T(t_N^*)e(i_1,\ldots,i_{N-1})  =
q^{\sum\limits_{k=1}^{N-1}i_k}e(i_1,\ldots,i_{N-1}).
\end{gather*}

Def\/ine the elements $\{x_j\}_{j=1,\ldots,N}$ as follows:
\begin{gather}\label{x_j}
x_j\overset{\mathrm{def}}{=}
\begin{cases}
\sum\limits_{k=j}^Nt_kt_k^*, & j>n,
\\ -\sum\limits_{k=j}^nt_kt_k^*+\sum\limits_{k=n+1}^Nt_kt_k^*, & j\le n.
\end{cases}
\end{gather}
Obviously, $x_1=1$, $x_ix_j=x_jx_i$,
\begin{gather}\label{t_jx_k}
t_jx_k=
\begin{cases}
q^2x_kt_j, & j<k,
\\ x_kt_j, & j\ge k.
\end{cases}
\end{gather}

The vectors $e(i_1,\ldots,i_{N-1})$ are joint eigenvectors of the operators $T(x_j)$, $j=1,2,\ldots,N$:
\begin{gather*}
T(x_j)e(i_1,\ldots,i_{N-1})= q^{2\sum\limits_{k=1}^{j-1}i_k}e(i_1,\ldots,i_{N-1}).
\end{gather*}

The joint spectrum of the pairwise commuting operators $T(x_j)$, $j=1,2,\ldots,N$, is
\begin{gather*}
\mathfrak{M}=\big\{(x_1,\ldots,x_N)\in\mathbb{R}^N \;| \; x_i/x_j\in q^{2\mathbb{Z}}\;\&\;1=x_1\le x_2\le\cdots\le
x_{n+1}, \\
\phantom{\mathfrak{M}=\big\{}{} \&\; x_{n+1}>x_{n+2}>\cdots>x_N>0\big\}.
\end{gather*}

The next proposition was proved in \cite{BS}.

\begin{proposition}
$T$ is a faithful representation of $\operatorname{Pol}(\mathscr{H}_{n,m})_q$.
\end{proposition}

Let us now introduce the notion of generalized functions on the quantum complex hyperbolic space
$\mathscr{H}_{n,m}$. Evidently, using the commutation relations, one can decompose every polynomial $f \in
\operatorname{Pol}(\widetilde{\mathscr{H}}_{n,m})_q$ as follows:
\begin{gather}\label{decomp_pol_Hnm_0}
f= \!\!\!\sum_{I=(i_1,\ldots,i_N), J=(j_1,\ldots,j_N) \in \mathbb Z_+^N}\!\!\! c_{IJ}t_1^{i_1}\cdots
t_n^{i_n}t_{n+1}^{*i_{n+1}}\cdots t_N^{*i_N}t_N^{j_N}\cdots t_{n+1}^{j_{n+1}}t_n^{*j_n}\cdots t_1^{*j_1},
\!\!\!\qquad c_{IJ} \in \mathbb C.\!\!\!
\end{gather}
Due to \eqref{x_j}, the latter can be reduced to the decomposition{\samepage
\begin{gather}\label{decomp_pol_Hnm}
f=\sum_{(i_1,\ldots,i_N,j_1,\ldots,j_N):\;i_kj_k=0} t_1^{i_1}\cdots t_n^{i_n}t_{n+1}^{*i_{n+1}}\cdots
t_N^{*i_N}f_{IJ}(x_2,\ldots,x_N)t_N^{j_N}\cdots t_{n+1}^{j_{n+1}}t_n^{*j_n}\cdots t_1^{*j_1},
\end{gather}
where $f_{IJ}(x_2,\ldots,x_N)$ are polynomials.}

One can equip $\operatorname{Pol}(\widetilde{\mathscr{H}}_{n,m})_q$ with the weakest topology such that the
functionals
\[
l_{(i_1,\ldots,i_N;j_1,\ldots,j_N)}(f)=(T(f)e(i_1,\ldots,i_N),e(j_1,\ldots,j_N))
\]
are continuous. The completion of $\operatorname{Pol}(\widetilde{\mathscr{H}}_{n,m})_q$ w.r.t.\ this topology
will be considered as the space of generalized functions on the quantum $\widehat{\mathscr{H}}_{n,m}$ and
denoted by $\mathscr{D}(\widetilde{\mathscr{H}}_{n,m})_q'$. Naturally, one can extend $T$ to a representation
of $\mathscr{D}(\widetilde{\mathscr{H}}_{n,m})_q'$ by continuity. Now \eqref{decomp_pol_Hnm} allows one to
identify $\mathscr{D}(\widetilde{\mathscr{H}}_{n,m})_q'$ with the space of formal series
\[
f=\sum_{(i_1,\ldots,i_N,j_1,\ldots,j_N):\;i_kj_k=0} t_1^{i_1}\cdots
t_n^{i_n}t_{n+1}^{*i_{n+1}}\cdots t_N^{*i_N}f_{IJ}(x_2,\ldots,x_N)t_N^{j_N}\cdots
t_{n+1}^{j_{n+1}}t_n^{*j_n}\cdots t_1^{*j_1},
\]
where $f_{IJ}(x_2,\ldots,x_N)$ are functions on $\mathfrak{M}$. The topology on this space of formal series
is the topology of pointwise convergence of the functions $f_{IJ}$.

Denote by $f_0$ the following function
\begin{gather}\label{f_0(x_n+1)}
f_0=f_0(x_{n+1})=
\begin{cases}
1, & x_{n+1}=1,
\\ 0, & x_{n+1}\in q^{-2\mathbb{N}}.
\end{cases}
\end{gather}
(Recall that $\mathrm{spec}\,x_{n+1}=q^{-2\mathbb{Z}_+}$.) Thus $f_0$ is a $q$-analog of the characteristic
function of the submanifold
\[
\big\{ (t_1,\ldots,t_N)\in\mathbb{C}^N\; \big|\;t_1=t_2=\cdots=t_n=0 \big\}\cap\mathscr{H}_{n,m}.
\]

Introduce now a $*$-algebra $\operatorname{Fun}(\widetilde{\mathscr{H}}_{n,m})_q \subset
\mathscr{D}(\widetilde{\mathscr{H}}_{n,m})_q'$ generated by
$\operatorname{Pol}(\widetilde{\mathscr{H}}_{n,m})_q$ and $f_0$. Easy computations from~\eqref{t_jx_k} show
that $f_0$ satisf\/ies the following relations:
\begin{gather*}
 t_j^*f_0=f_0t_j=0,\qquad j\le n,
\\  x_{n+1}f_0=f_0x_{n+1}=f_0,
\\  f_0^2=f_0^*=f_0,
\\  t_jf_0=f_0t_j,\qquad t_j^*f_0=f_0t_j^*,\qquad j\ge n+1.
\end{gather*}

The relation $I_\varphi f_0=f_0$ allows one to extend the $*$-automorphism $I_\varphi$ \eqref{I_phi} of the
algebra $\operatorname{Pol}(\widetilde{\mathscr{H}}_{n,m})_q$ to the $*$-automorphism of
$\operatorname{Fun}(\widetilde{\mathscr{H}}_{n,m})_q$. Let
\[
\operatorname{Fun}(\mathscr{H}_{n,m})_q\overset{\mathrm{def}}{=} \left\{ f\in
\operatorname{Fun}(\widetilde{\mathscr{H}}_{n,m})_q\;\big|\; I_\varphi f=f\right\}.
\]
Obviously, there exists a unique extension of the $*$-representation $T$ to a $*$-representation of the
$*$-algebra $\operatorname{Fun}(\mathscr{H}_{n,m})_q$ such that $T(f_0)$ is the orthogonal projection of
$\mathscr{H}$ onto the linear span of vectors $\{e(\underbrace{0,\ldots,0}_n,i_{n+1},\ldots,i_{N-1})|\:
i_{n+1},\ldots,i_{N-1}\in\mathbb{N}\}$.

Let $\mathscr{D}(\mathscr{H}_{n,m})_{q,\mathfrak{k}}$ be the two-sided ideal of
$\operatorname{Fun}(\mathscr{H}_{n,m})_q$ generated by $f_0$. We call this ideal the algebra of f\/inite
functions on the quantum hyperbolic space. It is a quantum analog for the algebra of $U \mathfrak{k}=U
\mathfrak{s(gl}_n \times \mathfrak{gl}_m)$-f\/inite smooth functions on $\mathscr{H}_{n,m}$ with compact
support.

\begin{remark}
Let us explain the adjective `f\/inite'. If $f$ is a f\/inite function, $T(f)$ is an operator with only a f\/inite
number of nonzero entries. However, we do not consider all possible f\/inite functions (and, therefore, all
operators with f\/inite number of nonzero entries) but only $U_q \mathfrak{k}$-f\/inite ones, cf.~\cite{Faraut}.
\end{remark}

It was proved in \cite{BS} that

\begin{theorem}
The representation $T$ of $\mathscr{D}(\mathscr{H}_{n,m})_{q,\mathfrak{k}}$ is faithful.
\end{theorem}

\begin{remark}\label{remark2}
Let $f(x_{n+1})$ be a polynomial. Then it follows from \eqref{x_j}, \eqref{t_jx_k} that
\begin{gather}\label{t_ift_^*}
\sum_{i=1}^nt_if(x_{n+1})t_i^*=f\big(q^2x_{n+1}\big)\sum_{i=1}^nt_it_i^*= f\big(q^2x_{n+1}\big)(x_{n+1}-1).
\end{gather}
This computation, together with \eqref{f_0(x_n+1)}, allows one to consider the element
$f_1=\sum\limits_{i=1}^nt_if_0t_i^*$ as a~function of $x_{n+1}$ such that
\[
f_1(x_{n+1})=
\begin{cases}
q^{-2}-1, & x_{n+1}=q^{-2},
\\ 0, & x_{n+1}=1\text{ \ or \ }x_{n+1}\in q^{-2\mathbb{N}-2}.
\end{cases}
\]

Thus a multiple application of \eqref{t_ift_^*} leads to the following claim:
$\mathscr{D}(\mathscr{H}_{n,m})_{q,\mathfrak{k}}$ contains {\it all} f\/inite functions of $x_{n+1}$ (i.e.,
such functions $f$ that $f(q^{-n})=0$ for all but f\/initely many $n\in\mathbb{N}$).
\end{remark}

The $U_{q}\mathfrak{su}_{n,m}$-module algebra structure is established on
$\mathscr{D}(\mathscr{H}_{n,m})_{q,\mathfrak{k}}$ by the following:
\begin{gather*}
E_nf_0 = -\frac{q^{-1/2}}{q^{-2}-1}t_nf_0t_{n+1}^*, \qquad F_nf_0 = -\frac{q^{3/2}}{q^{-2}-1}t_{n+1}f_0t_n^*,
\qquad K_nf_0 = f_0,
\\
E_jf_0 = F_jf_0=(K_j-1)f_0=0,\qquad j\ne n.
\end{gather*}

Now we present an explicit formula for a positive invariant integral on the space of f\/inite functions
$\mathscr{D}(\mathscr{H}_{n,m})_{q,\mathfrak{k}}$ and thereby establish its existence.

Let $\nu_q:\mathscr{D}(\mathscr{H}_{n,m})_{q,\mathfrak{k}} \to\mathbb{C}$ be a linear functional def\/ined by
\[
\nu_q(f)=\operatorname{Tr}(T(f)\cdot Q)=\int_{\mathscr{H}_{n,m}} fd\nu_q,
\]
where $Q:\mathscr{H}\to\mathscr{H}$ is the linear operator given on the basis elements
$e(i_1,\ldots,i_{N-1})$ by
\[
Qe(i_1,\ldots,i_{N-1})= \mathrm{const_1} \, q^{2\sum\limits_{j=1}^{N-1}(N-j)i_j}e(i_1,\ldots,i_{N-1}),
\qquad \mathrm{const_1}>0.
\]

\begin{theorem}[\protect{\cite{BS}}]
The functional $\nu_q$ is well defined, positive, and $U_{q}\mathfrak{su}_{n,m}$-invariant.
\end{theorem}

One has to normalize this integral in a some way. In~\cite{BS} we put $\int_{\mathscr{H}_{n,m}} f_0
d\nu_q=1$, so the constant in the previous theorem equals
\[
\mathrm{const_1}= \prod\limits_{j=1}^{m-1}\big(q^{-2j}-1\big).
\]

\section[The subalgebra $\mathscr{D}(\mathscr{H}_{n,m})_{q,\mathfrak{k}}^{U_q \mathfrak{k}}$]{The subalgebra $\boldsymbol{\mathscr{D}(\mathscr{H}_{n,m})_{q,\mathfrak{k}}^{U_q \mathfrak{k}}}$}\label{sec_4}

In this section we restrict ourselves to subalgebras of $U_q \mathfrak k$-invariant elements. It is well
known that $\mathscr{H}_{n,m}$ is a pseudo-Hermitian symmetric space of rank~1~\cite{Molch2}. The following
proposition is a~natural quantum analog for this fact.

\begin{proposition}
$\mathscr{D}(\mathscr{H}_{n,m})_{q,\mathfrak{k}}^{U_q \mathfrak{k}}=\{f(x_{n+1}) \in
\mathscr{D}(\mathscr{H}_{n,m})_{q,\mathfrak{k}}\}$.
\end{proposition}
Since $U_q \mathfrak{k}=\mathbb C[K_n,K_n^{-1}] \otimes U_q \mathfrak{sl}_n \otimes U_q \mathfrak{sl}_m$, the
proof of this proposition follows from the next statement on $U_q \mathfrak{sl}_n \times U_q
\mathfrak{sl}_m$-isotypic components in $\mathscr{D}(\mathscr{H}_{n,m})_{q,\mathfrak{k}}$.
\begin{proposition}
$U_q \mathfrak{sl}_n \times U_q \mathfrak{sl}_m$-isotypic components of
$\mathscr{D}(\mathscr{H}_{n,m})_{q,\mathfrak{k}}$ correspond to the modules $L^{(n)}(a\varpi_1+d\varpi_{n-1})
\boxtimes L^{(m)}(c\varpi_1 +b \varpi_{m-1})$ with $a,b,c,d \in \mathbb Z_+$, $a+c=b+d$ $($with infinite
multiplicity$)$. The highest weight subspace is spanned by the vectors
$t_1^at_N^{*b}\varphi(x_{n+1})t_{n+1}^ct_n^{*d}$, where $\varphi(x_{n+1})$ is a finite function.
\end{proposition}

\begin{remark}
The sign $\boxtimes$ ref\/lects the fact that the multipliers are modules of dif\/ferent algebras $U_q
\mathfrak{sl}_n$ and $U_q \mathfrak{sl}_m$.
\end{remark}

\begin{proof} Let us describe $U_q \mathfrak{sl}_n \times U_q \mathfrak{sl}_m$-highest weight vectors in
$\mathscr{D}(\mathscr{H}_{n,m})_{q,\mathfrak{k}}$. One can decompose every f\/inite function $f \in
\mathscr{D}(\mathscr{H}_{n,m})_{q,\mathfrak{k}}$ in the following way
\begin{gather*}
f= \sum_{\substack{I=(i_1,\ldots,i_N), J=(j_1,\ldots,j_N) \in \mathbb Z_+^N \\ i_1+\cdots +
i_n+j_{n+1}+\cdots +j_N\\ = i_{n+1}+\cdots+i_N+j_1+\cdots +j_n}} \!\!\!\!\!c_{IJ}t_1^{i_1}\cdots
t_n^{i_n}t_{n+1}^{*i_{n+1}}\cdots t_N^{*i_N}f_0t_N^{j_N}\cdots t_{n+1}^{j_{n+1}}t_n^{*j_n}\cdots t_1^{*j_1},
\qquad c_{IJ} \in \mathbb C.
\end{gather*}
(just by applying the commutation relations, cf.~\eqref{decomp_pol_Hnm_0}). Let $L^{(n)}(\lambda)$ be the
f\/inite dimensional $U_q \mathfrak{sl}_n$-module with highest weight $\lambda$ and $\varpi_j$ are fundamental
weights of $\mathfrak{sl}_n$. By standard arguments,
\begin{gather*}
 {\rm l.s.}\big\{t_1^{l_1}t_2^{l_2}\cdots t_n^{l_n}\,|\,l_1+\cdots+l_n=a\big\}=L^{(n)}(a\varpi_1),
\\  {\rm l.s.}\big\{t_1^{*l_1}t_2^{*l_2}\cdots t_n^{*l_n}\,|\,l_1+\cdots+l_n=a\big\}=L^{(n)}(a\varpi_{n-1}),
\\  {\rm l.s.}\big\{t_{n+1}^{l_{n+1}}t_{n+2}^{l_{n+2}}\cdots t_N^{l_N}\,|\, l_{n+1}+\cdots+l_N=a\big\}=L^{(m)}(a\varpi_1),
\\  {\rm l.s.}\big\{t_{n+1}^{*l_{n+1}}t_{n+2}^{*l_{n+2}}\cdots
t_N^{*l_N}\,|\,l_{n+1}+\cdots+l_N=a\big\}=L^{(m)}(a\varpi_{m-1}).
\end{gather*}
Thus we have an epimorphism
\begin{gather*}
\bigoplus_{a+c=b+d} L^{(n)}(a\varpi_1) \otimes  L^{(m)}(d\varpi_{m-1}) \otimes L^{(m)}(c\varpi_1) \otimes
L^{(n)}(b\varpi_{n-1}) \rightarrow \mathscr{D}(\mathscr{H}_{n,m})_{q,\mathfrak{k}},
\\
t_1^{i_1}\cdots t_n^{i_n} \boxtimes t_{n+1}^{*i_{n+1}}\cdots t_N^{*i_N} \otimes t_N^{j_N}\cdots
t_{n+1}^{j_{n+1}} \boxtimes t_n^{*j_n}\cdots t_1^{*j_1}\\
\qquad{} \mapsto t_1^{i_1}\cdots t_n^{i_n}
t_{n+1}^{*i_{n+1}}\cdots t_N^{*i_N}f_0t_N^{j_N}\cdots t_{n+1}^{j_{n+1}}t_n^{*j_n}\cdots t_1^{*j_1}.
\end{gather*}
Recall that, as in the classical case, see \cite{VO},
\[
L^{(n)}(a\varpi_1) \otimes L^{(n)}(b\varpi_{n-1}) = \oplus_{i=0}^{\max(a,b)}
L^{(n)}((a-i)\varpi_1+(b-i)\varpi_{n-1})
\]
in the category of $U_q \mathfrak{sl}_n$-modules. By explicit calculations one can show that the $U_q
\mathfrak{sl}_n$-highest weight vector in $L^{(n)}((a-i)\varpi_1+(b-i)\varpi_{n-1}) \subset
L^{(n)}(a\varpi_1) \otimes L^{(n)}(b\varpi_{n-1})$ has the form
\[
t_1^{a-i}  \left(\sum_{j=1}^nt_j
\otimes t_j^*\right)^{i}   t_n^{*(b-i)}.
\]

Now a routine application of the commutation relations (similar to the ones in Remark~\ref{remark2}) allows one to reduce
every highest weight vector to a linear span of the vectors $t_1^at_N^{*b}\varphi(x_{n+1})t_{n+1}^ct_n^{*d}$,
where $\varphi(x_{n+1})$ is a f\/inite function.
\end{proof}

Now let us obtain an explicit form of the restriction of $\nu_q$ to the space
$\mathscr{D}(\mathscr{H}_{n,m})_{q,\mathfrak{k}}^{U_q \mathfrak{k}}$ of $U_q \mathfrak{k}$-invariant elements
of $\mathscr{D}(\mathscr{H}_{n,m})_{q,\mathfrak{k}}$.

Recall the standard notation $(a;q)_k=(1-a)\cdots (1-aq^{k-1})$, $(a;q)_0=1$,
\[
\int_{1}^\infty f(x) d_{q^{-2}}x=\big(q^{-2}-1\big)\sum_{k=0}^\infty f\big(q^{-2k}\big)q^{-2k}.
\]

\begin{proposition}\label{rad_int}
For any function $f(x_{n+1}) \in \mathscr{D}(\mathscr{H}_{n,m})_{q,\mathfrak{k}}$ one has
\[
\int_{\mathscr{H}_{n,m}} f(x_{n+1}) d\nu_q= \int_1^\infty f(x)\rho(x) d_{q^{-2}}x,
\]
where
\begin{gather}\label{rho_x}
\rho(x)=\mathrm{const_2}\,   x^{m-1}\big(q^{-2}x-1\big)\big(q^{-4}x-1\big)\cdots \big(q^{-2(n-1)}x-1\big),
\end{gather}
and $\mathrm{const_2}=\frac{1}{q^{-2}-1}\prod \limits_{j=1}^{n-1}\frac{1}{q^{-2j}-1}.$
\end{proposition}

\begin{proof} By explicit calculations,
\begin{gather*}
\int_{\mathscr{H}_{n,m}} f(x_{n+1}) d\nu_q= \mathrm{const_1}
\sum_{\substack{i_1\ldots,i_n \in -\mathbb{Z}_+,\\
i_{n+1},\ldots,i_{N-1} \in \mathbb{N}}} f\big(q^{2i_1+\cdots+2i_n}\big)
q^{2(N-1)i_1+2(N-2)i_2+\cdots+2i_{N-1}}
\\
\qquad{} = \mathrm{const_1}    \sum \limits_{i_{n+1},\ldots,i_{N-1}=1}^\infty \left(\sum_{\substack{i_1, \ldots, i_n
\in -\mathbb Z_+, \\ i_1+\cdots+i_n=i}} f(q^{2i})
q^{2(N-1)i_1+\cdots+2mi_n}\right) q^{2(m-1)i_{n+1}+\cdots+2i_{N-1}}
\\
\qquad{}= \mathrm{const_1} \sum_{i \in -\mathbb Z_+} \frac{q^{m(m-1)}}{(q^2;q^2)_{m-1}} \left(\sum_{\substack{i_1, \ldots, i_{n-1}
\in -\mathbb Z_+, \\ i_1+\cdots+i_{n-1} \geq i}}
q^{2(n-1)i_1+\cdots+2i_{n-1}}\right)f(q^{2i})q^{2mi}
\\
\qquad{}= \mathrm{const_1} \frac{q^{m(m-1)}}{(q^2;q^2)_{m-1}} \sum_{i \in -\mathbb Z_+} f(q^{2i})q^{2mi}  \left(
\sum_{\substack{i_1, \ldots, i_{n-1} \in \mathbb Z_+, \\ i_1+\cdots+i_{n-1} \leq -i}}
q^{-2(n-1)i_1-\cdots-2i_{n-1}}\right).
\end{gather*}
Let us verify that for every $j,k \in \mathbb Z_+$
\[
\sum_{\substack{i_1, \ldots, i_k \in \mathbb Z_+, \\ i_1+\cdots+i_k
\leq j}}
q^{-2(k-1)i_1-\cdots-2i_k}=\frac{(q^{-2};q^{-2})_{j+k}}{(q^{-2};q^{-2})_j(q^{-2};q^{-2})_k}.
\]
Denote the l.h.s.\ of the previous equation by $\Psi(j,k)$. One can verify the recurrence relation for the
$q$-Pascal triangle
\[
\Psi(j,k)=q^{-2k}\Psi(j-1,k)+\Psi(j,k-1)
\]
and the boundary values
\[
\Psi(0,k)=1, \qquad \Psi(j,1)=\frac{1-q^{-2(j+1)}}{1-q^{-2}},
\]
explicitly. Thus $\Psi(j,k)$ are the corresponding entries of the $q$-Pascal triangle. Now we can complete our
calculations and obtain that
\begin{gather*}
\int_{\mathscr{H}_{n,m}} f(x_{n+1}) d\nu_q= \mathrm{const_1}
\frac{q^{m(m-1)}}{(q^2;q^2)_{m-1}} \sum_{i \in -\mathbb Z_+} f\big(q^{2i}\big)
q^{2mi}
\frac{(q^{-2};q^{-2})_{n-1-i}}{(q^{-2};q^{-2})_{-i}(q^{-2};q^{-2})_{n-1}}
\\
\qquad {} = \mathrm{const_1}   \frac{q^{m(m-1)}}{(q^2;q^2)_{m-1}} \sum_{i
\in -\mathbb Z_+} f\big(q^{2i}\big) q^{2mi}
\frac{(q^{2i-2};q^{-2})_{n-1}}{(q^{-2};q^{-2})_{n-1}}=\int_1^\infty
f(x)\rho(x) d_{q^{-2}}x,
\end{gather*}
where $\rho(x)=\mathrm{const_2} \,  x^{m-1}(q^{-2}x-1)(q^{-4}x-1)\cdots (q^{-2(n-1)}x-1),$ and
\begin{gather*}
\mathrm{const_2}=\frac{1}{q^{-2}-1} \mathrm{const_1}   \frac{q^{m(m-1)}}{(q^2;q^2)_{m-1}}
\frac{(-1)^{n-1}}{(q^{-2};q^{-2})_{n-1}} =
\frac{(-1)^{n-1}}{(q^{-2}-1)(q^{-2};q^{-2})_{n-1}}.  \tag*{\qed}
\end{gather*}
\renewcommand{\qed}{}
\end{proof}

Now one can complete $\mathscr{D}(\mathscr{H}_{n,m})_{q,\mathfrak{k}}^{U_q \mathfrak{k}}$ with respect to the
norm $||f||^2=\int_1^\infty f^*f \rho(x)d_{q^{-2}}x$. The resulting Hilbert space will be denoted by
$L^2\big(d\nu_q^{(0)}\big)$.

\section[Covariant first order differential calculi over $\operatorname{Pol}(\widetilde{\mathscr{H}}_{n,m})_q$]{Covariant f\/irst order dif\/ferential calculi over $\boldsymbol{\operatorname{Pol}(\widetilde{\mathscr{H}}_{n,m})_q}$}
\label{sec_5}

In this section we introduce   holomorphic and antiholomorphic covariant f\/irst order dif\/ferential calculi
over $\operatorname{Pol}(\widetilde{\mathscr{H}}_{n,m})_q$.

First of all, recall a general def\/inition of a covariant f\/irst order dif\/ferential calculus fol\-lowing~\cite{Kl-Sch}.

Let $F$ be a unital algebra. A f\/irst order dif\/ferential calculus over $F$ is a pair $(M,d)$ where $M$ is an
$F$-bimodule and $d:F \rightarrow M$ is a linear map such that
\begin{enumerate}\itemsep=0pt
\item for all $f_1,f_2 \in F$ one has $d(f_1 \cdot f_2)=df_1 \cdot f_2+f_1 \cdot df_2$;
\item $M$ is a linear span of the vectors $f_1 \cdot df_2 \cdot f_3$, where $f_1,f_2,f_3 \in F$.
\end{enumerate}

Now suppose that $A$ is a Hopf algebra and $F$ is an $A$-module algebra. A f\/irst order dif\/ferential calculus
$(M,d)$ over $F$ is called covariant if the following conditions hold:
\begin{enumerate}\itemsep=0pt
\item $M$ is an $A$-covariant $F$-bimodule, i.e.\ the action maps $F \otimes M \rightarrow M$, $M \otimes F \rightarrow M$
are morphisms of $A$-modules;
\item $d$ is a morphism of $A$-modules.
\end{enumerate}

Let $\mathcal{L}=L \circ \omega$, where $L$ is the canonical antirepresentation of $U_q \mathfrak{sl}_N$ in
$\mathbb C[SL_N]_q$ (or the representation of $U_q \mathfrak{sl}_N^{\rm op}$) and $\omega$ is an
antiautomorphism of $U_q \mathfrak{sl}_N$ def\/ined on the generators by $\omega(K_i)=K_i^{-1}$,
$\omega(E_i)=F_i$, $\omega(F_i)=E_i$. Thus $\mathcal{L}$ is a representation of $U_q \mathfrak{sl}_N$ in
$\mathbb C[SL_N]_q$.

Consider the action of $\mathcal{L}(E_1K_1^{-1/2})$ in $\mathbb C [SL_N]_q$ and its restriction to
$\operatorname{Pol}(\widetilde{\mathscr{H}}_{n,m})_q$. Put
\[\partial: \
\operatorname{Pol}(\widetilde{\mathscr{H}}_{n,m})_q \rightarrow \operatorname{Pol}(\widetilde{X})_q, \qquad \partial
t_i=\mathcal{L}\big(E_1K_1^{-1/2}\big)t_{1i}.
\]

Let $\Omega^{(1,0)}(\widetilde{\mathscr{H}}_{n,m})_q \subset \operatorname{Pol}(\widetilde{X})_q$ be the
$\operatorname{Pol}(\widetilde{\mathscr{H}}_{n,m})_q$-submodule generated by $\partial t_i$, $i=1,\dots,N$.
\begin{lemma}
$\Omega^{(1,0)}(\widetilde{\mathscr{H}}_{n,m})_q$ is a $U_q \mathfrak{sl}_N$-covariant first order differential calculus
over $\operatorname{Pol}(\widetilde{\mathscr{H}}_{n,m})_q$.
\end{lemma}
\begin{proof} One has to verify the Leibniz rule  which immediately follows from the formulas for the
comultiplication in $U_q \mathfrak{sl}_N$. Since left and right action of $U_q \mathfrak{sl}_N$ in $\mathbb C
[SL_N]_q$ commute, $\partial$ is a~morphism of $U_q \mathfrak{sl}_N$-modules, so the calculus is covariant.
\end{proof}

Now we def\/ine $\bar{\partial}: \operatorname{Pol}(\widetilde{\mathscr{H}}_{n,m})_q \rightarrow
\operatorname{Pol}(\widetilde{X})_q$ by the rule $\bar{\partial}f=(\partial f^*)^*$. Let
$\Omega^{(0,1)}(\widetilde{\mathscr{H}}_{n,m})_q \subset \operatorname{Pol}(\widetilde{X})_q$ be the
$\operatorname{Pol}(\widetilde{\mathscr{H}}_{n,m})_q$-submodule generated by $\bar{\partial} t_i$, $i=1,\dots,N$.
\begin{lemma}
$\Omega^{(0,1)}(\widetilde{\mathscr{H}}_{n,m})_q$ is a $U_q \mathfrak{sl}_N$-covariant first order differential calculus
over $\operatorname{Pol}(\widetilde{\mathscr{H}}_{n,m})_q$.
\end{lemma}

This lemma can be proved similarly to the previous one.
\begin{remark}
The introduced f\/irst order dif\/ferential calculi can be obtained in another
way. One should start from the canonical Wess--Zumino calculi on a quantum
complex space introduced in~\cite{RTF} and then turn to a localization of
the corresponding algebras of functions.
\end{remark}

Now we introduce a Hermitian pairing $\Omega^{(0,1)}(\widetilde{\mathscr{H}}_{n,m})_q \times
\Omega^{(0,1)}(\widetilde{\mathscr{H}}_{n,m})_q \rightarrow \operatorname{Pol}({\mathscr{H}}_{n,m})_q$.

Let $P: \operatorname{Pol}(\widetilde{X})_q \rightarrow \operatorname{Pol}({\mathscr{H}}_{n,m})_q$ be the
projection parallel to a sum of other $U_q \mathfrak{s}(\mathfrak{gl}_1 \times \mathfrak{gl}_{N-1})$-isotypic components of $\mathcal{L}$. Now we def\/ine
\[
(\theta_1,\theta_2)=P(\theta^*_2\theta_1), \qquad \theta_1,\theta_2 \in \Omega^{(0,1)}(\widetilde{\mathscr{H}}_{n,m})_q.
\]
By obvious commutativity of the left and right actions of $U_q \mathfrak{sl}_N$ we have

\begin{proposition}
The pairing $(\cdot,\cdot)$ is $U_q \mathfrak{sl}_N$-invariant.
\end{proposition}
\begin{lemma}
\[P(t_{2j}t^*_{2k})=q^{-2}\frac{1-q^2}{1-q^{2(N-1)}}(\varepsilon_{jk}-t_{1j}t^*_{1k}),\] where
\[
\varepsilon_{jk}=\begin{cases} q^{2(j-1)}, \quad & j=k, \  j \geq n+1, \\
-q^{2(j-1)}, \quad & j=k, \  j \leq n, \\
0, \quad & j \neq k.
\end{cases}
\]
\end{lemma}
\begin{proof} Since $U_q \mathfrak{s(gl}_1 \times \mathfrak{gl}_{N-1})=\mathbb C[K_1,K_1^{-1}] \otimes U_q
\mathfrak{sl}_{N-1}$, one can decompose the projection $P$ as follows: $P=P_0P_1$, where $P_1$ is a
projection to the subspace of $U_q \mathfrak{sl}_{N-1}$-invariant elements (w.r.t.\ the $\mathcal{L}$-action)
and $P_0$ is a projection to the subspace of elements that are preserved by the $\mathcal{L}(K_1)$-action.

Let $u_1,\dots,u_k$ be the standard basis in the $U_q \mathfrak{sl}_k$-module $L(\varpi_1)$, and
$v_1,\ldots,v_k$ the dual basis in the $U_q \mathfrak{sl}_k$-module $L(\varpi_{k-1})$, where $\varpi_1$ and
$\varpi_{k-1}$ are the fundamental weights. A standard argument on f\/inite dimensional $U_q
\mathfrak{sl}_k$-modules allows one to prove that $\sum\limits_{j=1}^k (-q)^{j-1}u_j \otimes v_j$ is a~$U_q
\mathfrak{sl}_k$-invariant element in the $U_q \mathfrak{sl}_k$-module $L(\varpi_1) \otimes L(\varpi_{k-1})$,
and the map
\begin{gather}\label{proj_inv}
u_i \otimes v_j \mapsto \begin{cases} \displaystyle \frac{1-q^2}{1-q^{2k}}(-q)^{i-1}\sum_{a=1}^k (-q)^{a-1} u_a \otimes v_a, \quad &  i=j,
\\ 0, \quad & i \neq j,
\end{cases}
\end{gather}
is a projection to the subspace of $U_q \mathfrak{sl}_k$-invariant elements parallel to other isotypic
components. By obvious reasons, for $j=1,\ldots,N-1$ the maps
\[
\phi_j: \ u_i \mapsto t_{i+1,j},\qquad \psi_j: \ v_i \mapsto \det \nolimits_q T_{i+1,j},
 \] admit extensions to morphisms of
$U_q \mathfrak{sl}_{N-1}$-modules
\[
\phi_j: \ L(\varpi_{N-2}) \rightarrow \mathbb C[SL_N]_q, \qquad \psi_j: \ L(\varpi_1) \rightarrow \mathbb
C[SL_N]_q
\]
(w.r.t.\ the $\mathcal{L}$-action). By def\/inition~\eqref{def_inv}, $t_{ij}^*=(-q)^{j-i}\det \nolimits_qT_{ij}$
for $i \leq m$ and $j \geq n+1$. Thus one can apply the map \eqref{proj_inv} to compute for $l \geq n+1$
\[
P_1(t_{2l}t_{2l}^*)=(-q)^{l-2}P_1(t_{2l} \det \nolimits_q T_{2l})
=(-q)^{l-2}\frac{1-q^2}{1-q^{2(N-1)}}\sum_{j=1}^{N-1}(-q)^{j-1}t_{j+1,l}\det \nolimits_q T_{j+1,l}.
\]
Since $\mathbb C[SL_N]_q$ is a Hopf algebra, one has $\sum\limits_{j=1}^{N} t_{lj}S(t_jk)=\delta_{lk}$ and
$S(t_{jk})=(-q)^{j-k}\det \nolimits_q T_{kj}$. Thus
\begin{gather*}
P_1(t_{2l}t_{2l}^*)=(-q)^{2(l-2)}\frac{1-q^2}{1-q^{2(N-1)}}\big( 1-(-q)^{1-l}t_{1l}\det \nolimits_q
T_{1l}\big)
\\
\phantom{P_1(t_{2l}t_{2l}^*)}{} = (-q)^{2(l-2)}\frac{1-q^2}{1-q^{2(N-1)}}\big( 1-(-q)^{-2(l-1)}t_{1l}t_{1l}^*\big).
\end{gather*}
The other cases can be verif\/ied in a similar way. \end{proof}

\begin{proposition}
\[
\big(\bar \partial x_{n+1}, \bar \partial x_{n+1}\big)= q^{2(n-1)} \frac{1-q^2}{1-q^{2(N-1)}}
x_{n+1}(1-q^{-2n}x_{n+1}).
\]
\end{proposition}

\begin{proof} An easy application of the previous lemma allows one to compute that
\[
\big(\bar \partial (t_{1k}^*t_{1l}), \bar \partial
(t_{1j}^*t_{1i})\big)=q^{-2}\frac{1-q^2}{1-q^{2(N-1)}}(\varepsilon_{jk}t_{1i}^*t_{1l}-(t_{1i}^*t_{1j})(t_{1k}^*t_{1l})).
\]
Thus for $j,k \geq n+1$ one has
\begin{gather*}
\big(\bar \partial (q^{-2j}t_{1j}^*t_{1j}), \bar \partial
(q^{-2k}t_{1k}^*t_{1k})\big)\\
\qquad{} =q^{-2}\frac{1-q^2}{1-q^{2(N-1)}}\big(\delta_{jk}q^{-2(j+1)}t_{1j}^*t_{1j}
-\big(q^{-2j}t_{1j}^*t_{1j}\big)\big(q^{-2k}t_{1k}^*t_{1k}\big)\big).
\end{gather*}
By easy computations, $q^{-2(n+1)}x_{n+1}=\sum\limits_{j=n+1}^N q^{-2j}t_{1j}^*t_{1j}$, which enables to prove the
claim.
\end{proof}

Let us f\/ix notation for $q$-dif\/ference operators:
\[
B_-: \ f(x) \mapsto \frac{f(q^{-2}x)-f(x)}{q^{-2}x-x}, \qquad B_+: \ f(x) \mapsto \frac{f(q^2x)-f(x)}{q^2x-x}.
\]
\begin{lemma}\label{lem_dbar}
\[
\big(\bar \partial f(x_{n+1}), \bar \partial g(x_{n+1})\big)= q^{2(n-1)} \frac{1-q^2}{1-q^{2(N-1)}} x_{n+1}
\big(1-q^{-2n}x_{n+1}\big) \overline{B_-g(x_{n+1})} B_-f(x_{n+1}).
\]
\end{lemma}
\begin{proof} By explicit calculations in $\mathbb C[SL_N]_q$ we obtain for $i \geq n+1$
\[
t_{1i}\bar \partial x_{n+1}= q^{-1} \bar \partial x_{n+1} t_{1i}, \qquad t_{1i}^* \bar \partial
x_{n+1}=q^{-1} \bar
\partial x_{n+1} t_{1i}^*.
\]
Let us verify the second identity. One has $\bar \partial x_{n+1}=\sum\limits_{j=n+1}^Nt_{1j}t_{2j}^*$, so
\begin{gather*}
t_{1i}^*\bar \partial x_{n+1}=t_{1i}^*  \left(\sum_{j=n+1}^{i-1}
t_{1j}t_{2j}^*+t_{1i}t_{2i}^*+\sum_{j=i+1}^Nt_{1j}t_{2j}^*\right) = q^{-1} \sum_{j=n+1}^{i-1}t_{1j}t_{2j}^*
  t_{1i}^*
\\
{}
+q^{-1}t_{1i}t_{2i}^*t_{1i}^*+\big(q-q^{-1}\big)\sum_{j>i}t_{1j}t_{1j}^*t_{2i}^*+q^{-1}\sum_{j=i+1}^N
t_{1j}t_{2j}^*t_{1i}^*+\big(q^{-1}-q\big)\sum_{j>i}t_{1j}t_{1j}^*t_{2i}^*.
\end{gather*}
So we have $x_{n+1}\bar \partial x_{n+1}=q^{-2} \bar
\partial x_{n+1} x_{n+1}$. Hence for every polynomial $f(x)$
\[
\bar \partial f(x_{n+1})=\bar \partial x_{n+1}\frac{f(q^{-2}x_{n+1})-f(x_{n+1})}{q^{-2}x_{n+1}-x_{n+1}}=\bar
\partial x_{n+1} B_-(f)(x_{n+1}),
\]
and  the claim follows from the previous proposition. \end{proof}

Using the above formal arguments, we
can extend $\bar \partial$ to
$\mathscr{D}(\mathscr{H}_{n,m})_{q,\mathfrak{k}}$.

\section[The Laplace-Beltrami operator and its radial part]{The Laplace--Beltrami operator and its radial part}\label{sec_6}

In this section we introduce a $U_q \mathfrak{sl}_N$-invariant operator $\square$ in the space of f\/inite
functions $\mathscr{D}(\mathscr{H}_{n,m})_{q, \mathfrak{k}}$. This operator will be considered as a quantum
analog for the invariant Laplace--Beltrami operator on the complex hyperbolic space $\mathscr{H}_{n,m}$. Also
we compute an explicit formula for the restriction $\square^{(0)}$ of $\square$ to the space
$\mathscr{D}(\mathscr{H}_{n,m})_{q,\mathfrak{k}}^{U_q \mathfrak{k}}$, the so called radial part of~$\square$.
The restriction appears to be a $q$-dif\/ference operator in variable $x=x_{n+1}$.

Now we def\/ine $\square$ by the formula
\[
\int_{\mathscr{H}_{n,m}}f_2^* (\square f_1) d\nu_q=\int _{\mathscr{H}_{n,m}} (\bar \partial f_1, \bar
\partial f_2) d\nu_q, \qquad f_1,f_2 \in \mathscr{D}(\mathscr{H}_{n,m})_{q,\mathfrak{k}}.
\]

\begin{proposition}
$\square$ is a self-adjoint $U_q \mathfrak{sl}_N$-invariant operator.
\end{proposition}

\begin{proof} The self-adjointness follows from the def\/inition. The $U_q \mathfrak{sl}_N$-invariance of the
form $(\cdot,\cdot)$ and the linear functional $\int_{\mathscr{H}_{n,m}} \cdot d\nu_q$ implies the $U_q
\mathfrak{sl}_N$-invariance of~$\square$.
\end{proof}

Let us f\/ind the restriction $\square^{(0)}$ of the operator $\square$ on the subspace
\[
\mathscr{D}(\mathscr{H}_{n,m})_{q,\mathfrak{k}}^{U_q \mathfrak{k}}=\big\{ f(x_{n+1}), \ \mathrm{supp} f \subset q^{-2 \mathbb{Z}_+},\
\sharp (\mathrm{supp} f) < \infty \big\},
\]
based on the equation
\begin{gather*}
\int_1^\infty \big(\square^{(0)} f\big) (x) \overline{f(x)} \rho(x) d_{q^{-2}} x \\
\qquad{} = q^{2(n-1)}
\frac{1-q^2}{1-q^{2(N-1)}} \int_1^\infty x \big(1-q^{-2n} x\big) \left| \frac{f(q^{-2}x)-f(x)}{q^{-2}x - x}
\right|^2 \rho(x) d_{q^{-2}} x
\end{gather*}
for every $f$ with f\/inite support.

Using the $q$-analog of the partial integration formulas one can prove that operators $B_-$ and $-q^2 B_+$
are formally dual. Exactly, for every functions $u(x)$, $v(x)$ with f\/inite support on $q^{2\mathbb{Z}}$ holds
\[ \int_0^\infty u(x) \frac{v(q^{-2}x)-v(x)}{(q^{-2}-1)x} d_{q^{-2}} x = -q^2 \int_0^\infty
\frac{u(x)-u(q^2x)}{(1-q^2)x} v(x) d_{q^{-2}} x,
\]
where
\[
\int_0^\infty f(x) d_{q^{-2}} x = \big(q^{-2} - 1\big)
\sum\limits_{k=-\infty}^\infty f\big(q^{-2k}\big) q^{-2k}.
\]

Thus
\[
\square^{(0)}: \ f(x) \mapsto \mathrm{const}(q, n, N) \rho(x)^{-1} B_+ x \big(q^{-2n} x - 1\big) \rho(x) B_- f(x),
\]
where
\[
\mathrm{const}(q, n, N) = q^{2n} \frac{1-q^2}{1-q^{2(N-1)}},
\]
and $\rho(x)$ is def\/ined in \eqref{rho_x}.

Lemma~\ref{lem_dbar} allows us to extend the Hermitian $U_q
\mathfrak{sl}_N$-invariant pairing to f\/irst order dif\/ferential forms with
coef\/f\/icients in $\mathscr{D}(\mathscr{H}_{n,m})_{q,\mathfrak{k}}$.

\begin{lemma}
For every function $f(x)$ with finite support on $q^{-2\mathbb{Z}_+}$
\begin{gather*}
\square^{(0)} f(x) = \frac{q^2}{(1-q^2)^2x} \big((x-1)q^{2m-2}f(q^2x) +
\big(q^{-2n}x-1\big)f\big(q^{-2}x\big)\\
\phantom{\square^{(0)} f(x) =}{}
+\big(1+q^{2m-2}-q^{2m-2}x-q^{-2n}x\big)f(x)\big).
\end{gather*}
\end{lemma}

\begin{proposition}\label{cont_spec}
The operator $\square^{(0)}$ is bounded.
\end{proposition}

\begin{proof} Consider the operator
\[
\widetilde{\square} : \ f(x) \mapsto x^{-(N-1)} \mathrm{const}(q, n, N) q^{-2n} B_+ x^{N+1} B_-f(x).
\]
It dif\/fers from $\square^{(0)}$ by a compact operator, so it is suf\/f\/icient to prove that
$\widetilde{\square}$ is bounded. The boundness of the latter operator can be proved by the direct
evaluation.
\end{proof}

By the previous proposition, one can extend $\square^{(0)}$ from
$\mathscr{D}(\mathscr{H}_{n,m})_{q,\mathfrak{k}}^{U_q \mathfrak{k}}$ to a bounded self-adjoint operator in
$L^2\big(d\nu_q^{(0)}\big)$. The latter extension will also be denoted by $\square^{(0)}$.

\section[Generalized eigenfunctions of $\square^{(0)}$ and Al-Salam-Chihara polynomials]{Generalized eigenfunctions of $\boldsymbol{\square^{(0)}}$\\ and Al-Salam--Chihara polynomials}\label{sec_7}

In this section we obtain the initial results on the bounded self-adjoint operator $\square^{(0)}$, namely we
obtain its formal eigenfunctions and eigenvalues explicitly. Note that explicit computations of the
asymptotics of these eigenfunctions, as in the classical case, allow us to consider a quantum analog of the
Harish-Chandra $c$-function (see Appendix~\ref{appendixA}).

By the direct computation we obtain the following lemma.

\begin{lemma}\label{c1.5.2} The function
\[
\Phi_l(x) = {_{3}\Phi_{2}}\left(
\begin{array}{c}
x,q^{-2l},q^{2(l+N-1)}\\ q^{2n},0
\end{array}
;q^2,q^2\right)
\]
in $\mathscr{D}(\mathscr{H}_{n,m})_q^\prime$ is a generalized eigenfunction for $\square^{(0)}$:
\[\square^{(0)} \Phi_l=\lambda(l)\Phi_l\] with the eigenvalue
\[
\lambda(l) = - q^{2-2n}\frac{(1-q^{-2l})(1-q^{2l+2(N-1)})}{(1-q^2)^2}.
\]
\end{lemma}

Recall the def\/inition of the Al-Salam--Chihara polynomials, following \cite{koekoek}:
\[
Q_k(z;a,b|q)=\frac{(ab;q)_k}{a^k}\ {_{3}\Phi_{2}}\left(
\begin{array}{c}
q^{-k},a e^{i\theta},ae^{-i\theta}\\ ab,0
\end{array}
;q,q\right), \qquad z=\cos\theta.
\]

We have the following well known orthogonality relations for the Al-Salam--Chihara polynomials:
\begin{gather}
\frac{1}{2\pi}\int_{-1}^1 Q_i(z;a,b|q)Q_j(z;a,b|q)\frac{w(z)}{\sqrt{1-z^2}}\,dz\nonumber\\
\qquad{}+\sum_{1<aq^k<a}w_k Q_i(z_k;a,b|q)Q_j(z_k;a,b|q)=\frac{\delta_{ij}}{(q^{i+1},abq^i;q)_\infty},\label{rel_orth_ASC}
\end{gather}
where
\[
w(z)=\frac{h(z,1)h(z,-1)h(z,q^{1/2})h(z,-q^{1/2})}{h(z,a)h(z,b)}, \!\!\qquad
h(z,a)=\big(ae^{i\theta},ae^{-i\theta};q\big)_{\infty},\!\! \qquad z=\cos\theta,
\]
$z_k=\frac{aq^k+aq^{-k}}{2}$, and
\[
w_k=\frac{(a^{-2};q)_\infty}{(q,ab,b/a;q)_\infty}\frac{(1-a^2q^{2k})(a^2,ab;q)_k}
{(1-a^2)(q,qa/b;q)_k}q^{-k^2}\left(\frac{1}{a^3b}\right)^k.
\]

Also, these polynomials satisfy the following recurrence relations:
\begin{gather}
zQ_i(z;a,b|q)=\frac12
Q_{i+1}(z;a,b|q)+\frac12q^i(a+b)Q_i(z;a,b|q)\nonumber\\
\phantom{zQ_i(z;a,b|q)=}{} +\frac12\big(1-q^i\big)\big(1-abq^{i-1}\big)Q_{i-1}(z;a,b|q).\label{rec_ASC}
\end{gather}

As one can see, the eigenfunctions $\Phi_l(x)$ are connected with the Al-Salam--Chihara polynomials:
\[
\Phi_l(q^{-2k}) = \frac{q^{k(N-1)}}{(q^{2n};q^2)_k}Q_k\big(z;q^{2n-N+1},q^{N-1}|q^2\big), \qquad
e^{i\theta}=q^{2l+N-1},\qquad z=\cos\theta.
\]

Let us denote the orthogonal measure for Al-Salam--Chihara polynomials by $d\sigma$.  Note that if $2n-N+1>0$,
the sum in the orthogonality relations vanishes and we have just the continuous measure.

\section[A spectral theorem for the radial part of $\square$]{A spectral theorem for the radial part of $\boldsymbol{\square}$}\label{sec_9}

In this section we obtain a spectral theorem for $\square^{(0)}$. As in the classical case~\cite[pp.~429--432]{Faraut}, the support of the Plancherel measure consists of continuous and discrete parts.
The continuous part corresponds to principal unitary series of $U_q \mathfrak{su}_{n,m}$-modules related to a
quantum analog of the cone
\[
\Xi_{n,m}=\big\{x \in \mathbb C^{N}\,|\, {-}x_1 \bar{x}_1-\cdots -x_n \bar{x}_n + x_{n+1}\bar{x}_{n+1}+\cdots
+x_N\bar{x}_N=0\big\}
\]
(this series is established in~\cite{BS}). The discrete part is supposed to be corresponded to a discrete
series of unitary $U_q \mathfrak{su}_{n,m}$-modules.

\begin{theorem}
The bounded self-adjoint linear operator $\square^{(0)}$ is unitary equivalent to the operator of
multiplication by independent variable in the Hilbert space $L^2(d\sigma)$. The unitary equivalence is given
by the operator
\begin{gather*}
U:\  L^2\big(d\nu_q^{(0)}\big)  \rightarrow L^2(d\sigma),
\\
U:\ f(x) \mapsto\hat{f}(\lambda)=\int_1^\infty f(x)\Phi_l(x)\rho(x)d_{q^{-2}}x,
\end{gather*}
where $\lambda=-q^{2-2n}\frac{(1-q^{-2l})(1-q^{2l+2N-2})}{(1-q^2)^2}$.
\end{theorem}

\begin{proof} Let us consider the f\/inite functions on $q^{-2\mathbb Z_+}$
\[
f_j(x)=
\begin{cases} 1, & x=q^{-2j},\\ 0, & \text{otherwise}, \end{cases}
\qquad j \in \mathbb Z_+.
\]
These functions form an orthogonal system in $\mathscr{D}(\mathscr{H}_{n,m})_{q,\mathfrak{k}}^{U_q \mathfrak{k}}$
and the completion of its linear span is $L^2\big(d\nu_q^{(0)}\big)$. By standard arguments~\cite{AG}, the bounded
self-adjoint linear operator $\square^{(0)}$ is unitary equivalent to the multiplication operator $f(\lambda)
\mapsto \lambda f(\lambda)$ in the Hilbert space $L^2(d\mu(\lambda))$ of square integrable functions with
respect to a certain measure $d\mu(\lambda)$ with compact support in $\mathbb R$. Let us f\/ind explicitly the
corresponding measure and the operator of unitary equivalence $U$. One can f\/ix the unitary equivalence
operator by the condition $Uf_0=1$.

By easy calculations, $\square^{(0)}$ transforms $f_j$ by the formula
\begin{gather*}
\square^{(0)}f_j=\frac{q^2}{(1-q^2)^2} \big( f_{j+1}\big(1-q^{2j+2}\big)q^{2m-2}+f_{j-1}\big(1-q^{2j+2n-2}\big)q^{-2n}\\
\phantom{\square^{(0)}f_j=}{} +f_j
\big(q^{2j}+q^{2j+2m-2}-q^{2m-2}-q^{-2n}\big)\big), \qquad j \in \mathbb Z_+,
\end{gather*}
and $||f_j||^2=q^{-2j(N-1)}\frac{(q^{2j+2};q^2)_{n-1}}{(q^2;q^2)_{n-1}}$ (here we naturally suppose
$f_{-1}=0$).

{\sloppy Thus the f\/inite functions $e_j=q^{j(N-1)}\sqrt{\frac{(q^{2};q^2)_{n-1}}{(q^{2j+2};q^2)_{n-1}}}f_j$ form an
orthonormal system in $L^2\big(d\nu_q^{(0)}\big)$ and $\square^{(0)}$ acts on them by the formula
\begin{gather*}
\square^{(0)}e_j=\frac{q^{m-n+1}}{(1-q^2)^2}\Big(
e_{j+1}\sqrt{\big(1-q^{2j+2}\big)\big(1-q^{2j+2n}\big)}+e_{j-1}\sqrt{\big(1-q^{2j}\big)\big(1-q^{2j+2n-2}\big)} \\
\phantom{\square^{(0)}e_j=}{}+e_j
\big(q^{2j+(N-1)}+q^{2j+2n-(N-1)}-q^{N-1}-q^{1-N}\big)\Big), \qquad j \in \mathbb Z_+.
\end{gather*}

}

Thus $P_j=Ue_j \in L^2(d\mu(\lambda))$, $j\in\mathbb{Z}_+$ form an orthonormal system of polynomials
\[
\int P_i(\lambda)P_j(\lambda)d\mu(\lambda)=\delta_{ij},\qquad i,j\in\mathbb{Z}_+,
\]
and one has{\samepage
\begin{gather}
\lambda P_j(\lambda)= \frac{q^{m-n+1}}{(1-q^2)^2}
\Big(P_{j+1}(\lambda)\sqrt{\big(1-q^{2j+2}\big)\big(1-q^{2j+2n}\big)}\nonumber\\
\phantom{\lambda P_j(\lambda)=}{}
 +P_{j-1} (\lambda) \sqrt{\big(1-q^{2j}\big)
\big(1-q^{2j+2n-2}\big)}\nonumber\\
\phantom{\lambda P_j(\lambda)=}{}
+P_j(\lambda) \big(q^{2j+(N-1)}+q^{2j+2n-(N-1)}-q^{N-1}-q^{1-N}\big)\Big), \qquad j \in \mathbb Z_+,\label{P_rec}
\\
\label{init_P0}
P_0(\lambda)=Ue_0=Uf_0=1.
\end{gather}
(we naturally suppose that $P_{-1}=0$).}

The orthogonal polynomials $P_j(\lambda)$, $j\in\mathbb{Z}_+$ are determined by~\eqref{P_rec} and
\eqref{init_P0}. Let us compare them with the corresponding recurrence relations~\eqref{rec_ASC} and the
initial data for the Al-Salam--Chihara polynomials $Q_j(z;q^{2n-(N-1)},q^{N-1}|q^2)$. One can observe that
\[
||Q_j||^2=\frac{1}{(q^{2j+2},q^{2j+2n};q^2)_\infty}=\frac{(q^2,q^{2n};q^2)_j}{(q^2;q^2)_{n-1}(q^{2n};q^2)^2_\infty}.
\]

The polynomials
\[
\widetilde{Q}_j\overset{\mathrm{def}}{=}\sqrt{\frac{1}{(q^2,q^{2n};q^2)_j}} \,Q_j
\]
satisfy the following recurrence relations:
\begin{gather*}
z\widetilde{Q}_j=\dfrac12\sqrt{\big(1-q^{2(j+1)}\big)\big(1-q^{2j+2n}\big)}\widetilde{Q}_{j+1}+
\dfrac12\sqrt{\big(1-q^{2j}\big)\big(1-q^{2j+2n-2}\big)}\widetilde{Q}_{j-1}\\
\phantom{z\widetilde{Q}_j=}{} +\dfrac12q^{2j}\big(q^{2n-N+1}+q^{N-1}\big)\widetilde{Q}_j, \qquad j \in \mathbb Z_+.
\end{gather*}

Thus we obtain that $P_j$ and $\widetilde{Q}_j$ are related by the change of variable
\[
\lambda=\frac{2q^{m-n+1}}{(1-q^2)^2}z-\frac{q^{m-n+1}}{(1-q^2)^2}\big(q^{N-1}+q^{-(N-1)}\big).
\]
So,
\[
P_j(\lambda)=\sqrt{\frac{1}{(q^2,q^{2n};q^2)_j}}  Q_j\big(z;q^{2n-N+1},q^{N-1}|q^2\big),
\]
where $z=\cos \theta=\frac{e^{i\theta}+e^{-i\theta}}{2}$ and
\begin{gather*}
\lambda=-\frac{q^{2-2n}}{(1-q^2)^2}\big(1-q^{N-1}e^{i\theta}\big)\big(1-q^{N-1}e^{-i\theta}\big),\\
Uf_j=q^{-j(N-1)}\sqrt{\frac{(q^{2j+2};q^2)_{n-1}}{(q^{2};q^2)_{n-1}}}\sqrt{\frac{1}{(q^2,q^{2n};q^2)_j}}
Q_j\big(z;q^{2n-N+1},q^{N-1}|q^2\big) \\
\phantom{Uf_j}{} =q^{-j(N-1)}\frac{1}{(q^2;q^2)_j} Q_j\big(z;q^{2n-N+1},q^{N-1}|q^2\big).
\end{gather*}
On the other hand,
\[
\int_1^\infty
f_j(x)\Phi_{l}(x)\rho(x)d_{q^{-2}}x=\frac{q^{-j(N-1)}}{(q^2;q^2)_j} Q_j\big(z;q^{2n-N+1},q^{N-1}|q^2\big),
\]
where $z=\frac12\big(q^{2l+N-1}+q^{-(2l+N-1)}\big)$.

Hence for every function $f(x)$ on $q^{-2\mathbb{Z}_+}$ with f\/inite support one has
\[
Uf=\int_1^\infty f(x)\Phi_{l}(x)\rho(x)d_{q^{-2}}x.
\]
Now the claim of the theorem follows from the orthogonality relations for the Al-Salam--Chihara polynomials~\eqref{rel_orth_ASC}. \end{proof}

\begin{remark}
There is an extensive literature on harmonic analysis related to a quantum $SU_q(1,1)$ (see the references
\cite{KSR,U,KV} from three research groups and references therein). Even the notion of quantum $SU_q(1,1)$
had some uncertainties back then (unlike to the case of quantum compact group $SU_q(2)$). A remarkable quantum ef\/fect lies in the fact that the
Plancherel type theorems for the quantum $SU_q(1,1)$ present decompositions with entries coming from the
principal unitary series and strange unitary series of representations (which also has no classical analog).

Later on L.~Vaksman established a new concept of quantum $SU_q(1,1)$ (and other noncompact real Lie groups).
He substituted the group $SU(1,1) \subset SL_2(\mathbb C)$ by its principal homogeneous space $w_0SU(1,1)$,
where $w_0=\begin{pmatrix}0, & -1 \\ 1, & 0\end{pmatrix}$. This enables him and his collaborators to prove
the Plancherel type theorems with decompositions that does not contain the strange series entries (see~\cite{SSV,BK}). Still a relation between these two approaches of quantization may require further
investigations.
\end{remark}

\appendix

\section[Appendix on the Harish-Chandra $c$-function and the Plancherel measure]{Appendix on the Harish-Chandra $\boldsymbol{c}$-function\\ and the Plancherel measure}\label{appendixA}

Let us introduce the notation of the $q$-analog of the Harish-Chandra $c$-function
\[
c(l)=\frac{\Gamma_{q^2}(n)\Gamma_{q^2}(2l+N-1)}{\Gamma_{q^2}(l+n)\Gamma_{q^2}(l+N-1)}=
\frac{(q^{2(l+N-1)};q^2)_\infty(q^{2(l+n)};q^2)_\infty}{(q^{2(2l+N-1)};q^2)_\infty(q^{2n};q^2)_\infty},
\]
where $\Gamma_q(x)=(1-q)^{1-x}\frac{(q;q)_\infty}{(q^x;q)_\infty}$ is a well known quantum analog for the
Gamma-function \cite{GR}.
\begin{lemma}\label{assymp}\qquad
\begin{enumerate}\itemsep=0pt
\item[$1.$] For $\mathrm{Re}\, l < -(N-1)/2$ one has $\Phi_l(x) \sim
c(-l-N+1)x^{-l-N+1}$ as $x \in q^{-2\mathbb Z_+}$ tends to infinity.
\item[$2.$] For $\mathrm{Re}\, l > -(N-1)/2$
one has $\Phi_l(x) \sim c(l)x^l$ as $x \in q^{-2\mathbb Z_+}$ tends to infinity.
\end{enumerate}
\end{lemma}

\begin{proof} Let us compute an asymptotic behavior of $\Phi_l(x)$ for $\mathrm{Re} \,l > -(N-1)/2$, the other
case can be managed similarly. For $x=q^{-2k}$, $k \rightarrow \infty$, one has
\begin{gather*}
{_{3}\Phi_{2}}\left( \begin{array}{c} q^{-2k},q^{-2l},q^{2(l+N-1)}\\
q^{2n},0
\end{array} ;q^2,q^2\right)= \frac{(q^{2-2l-2n-2k};q^2)_k}{(q^{2-2n-2k};q^2)_k} {_{2}\Phi_{1}}\left(
\begin{array}{c} q^{-2k},q^{-2l}\\ q^{2-2l-2n-2k}
\end{array} ;q^2,q^{2l+2m}\right)
\\
\hphantom{{_{3}\Phi_{2}}\left( \begin{array}{c} q^{-2k},q^{-2l},q^{2(l+N-1)}\\
q^{2n},0
\end{array} ;q^2,q^2\right)}{}
\sim q^{-2lk}
\frac{(q^{2l+2n};q^2)_\infty}{(q^{2n};q^2)_\infty} {_{1}\Phi_{0}}\left(
\begin{array}{c} q^{-2l}\\ - \end{array} ;q^2,q^{4l+2(N-1)}\right) \\
\hphantom{{_{3}\Phi_{2}}\left( \begin{array}{c} q^{-2k},q^{-2l},q^{2(l+N-1)}\\
q^{2n},0
\end{array} ;q^2,q^2\right)}{}
=q^{-2lk}
\frac{(q^{2l+2n};q^2)_\infty(q^{2l+2(N-1)};q^2)_{\infty}}{(q^{2n};q^2)_\infty(q^{4l+2(N-1)};q^2)_\infty}.
  \tag*{\qed}
\end{gather*}
\renewcommand{\qed}{}
\end{proof}

Let us return to the Plancherel measure of $\square^{(0)}$, namely for its absolutely continuous part. Since
$e^{i\theta}=q^{2l+N-1}$ and $z=\cos \theta$, one can observe by easy calculations that
\begin{gather*}
w(z)=w\big(z,q^{2n-N+1},q^{N-1}|q^2\big)=\frac{(q^{4l+2(N-1)},q^{-4l-2(N-1)};q^2)_\infty}{(q^{2l+2n},q^{2n-2l-2(N-1)},
q^{2l+2(N-1)},q^{-2l};q^2)_\infty}
\\
\hphantom{w(z)=w\big(z,q^{2n-N+1},q^{N-1}|q^2\big)}{}
=\frac{1}{c(l)c(-l-(N-1))}\frac{1}{(q^{2n};q^2)^2_\infty}.
\end{gather*}
By this calculations we establish a relation between continuous part of the Plancherel measure and the
Harish-Chandra $c$-function in the quantum case. Of course, this interplay in the classical ($q=1$) situation
is well-known.

\subsection*{Acknowledgements}

This project started out as joint work with L.~Vaksman and D.~Shklyarov. We are grateful to both of them for
helpful discussions and drafts with preliminary def\/initions and computations. Also we are grateful for
referees for their comments that help to improve and simplify our exposition.

\pdfbookmark[1]{References}{ref}
\LastPageEnding

\end{document}